\documentclass[11pt,leqno]{article}
\usepackage{amsmath,amssymb,float,tikz}
\usepackage{graphicx}

\def\nexto{\kern -0.54em}

\newcommand{\Z}{{\mathbb Z}}

\newcommand{\R}{{\mathbb R}}

\newcommand{\ttt}{|\hspace{-0.25mm}|\hspace{-0.25mm}|}
\renewcommand{\Re}{\mathop{\rm Re}\nolimits}
\renewcommand{\Im}{\mathop{\rm Im}\nolimits}

\begin{document}

\begin{center}
  {\Large\bf Semi-continuity of Oseledets flags\\ and Pesin sets with exponentially small tails} 
\end{center}

\begin{center}
{ Luchezar Stoyanov}
\end{center}


\begin{center}
{\it University of Western Australia,Crawley 6009 WA, Australia\\(e-mail:luchezar.stoyanov@uwa.edu.au)}
\end{center}

\noindent
{\it Abstract.}
Let $f$ be an invertible  transitive subshift of finite type over a bilateral symbol space $X$, let $\mu$ be a Gibbs measure 
for $f$ determined by a H\"older continuous potential on $X$, and let $A$ be an invertible  continuous linear cocycle over $f$ acting on 
a continuous $\R^d$-bundle $E$ over $X$ with Lyapunov exponents $\lambda_k < \lambda_{k-1} < \ldots < \lambda_1$
such that $A^{-1}$ is continuous as well.
We prove that if the Oseledets flags $F_j(x) = E_j(x) \oplus E_{j-1}(x) \oplus \ldots \oplus E_1(x)$
depend upper semi-continuously on $x \in X$,
then there exists a Pesin set with exponentially small tails for $\mu$. 

\medskip 


\medskip


\def\nexto{\kern -0.54em}

\def\T{{\bf T}}
\def\S{{\mathbb S}}
\def\diam{\mbox{\rm diam}}
\def\rr{{\mathcal R}}
\def\mt{{\Lambda}}
\def\e{\emptyset}
\def\sd{\S^{d-1}}
\def\so{\S^1}
\def\dQ{\partial Q}
\def\dk{\partial K}
\def\endofproof{{\rule{6pt}{6pt}}}
\def\ts{\tilde{\sigma}}
\def\tr{\tilde{r}}
\def\sigt{\ts^{r}}
\def\di{\displaystyle}
\def\dist{\mbox{\rm dist}}
\def\tU{\widetilde{U}}
\def\hU{\widehat{U}}
\def\tS{\tilde{S}}
\def\tP{\widetilde{\Pi}}
\def\sa+{\Sigma_A^+}
\def\u-{\overline{u}}
\def\du{\frac{\partial}{\partial u}}
\def\dv{\frac{\partial}{\partial v}}
\def\dt{\frac{d}{d t}}
\def\dx{\frac{\partial}{\partial x}}
\def\z-{\overline{z}}
\def\te{\tilde{e}}
\def\tv{\tilde{v}}
\def\tu{\tilde{u}}
\def\tw{\tilde{w}}
\def\con{\mbox{\rm const }}
\def\Box{\spadesuit}
\def\nn{{\mathcal N}}
\def\mm{{\mathcal M}}
\def\kk{{\mathcal K}}
\def\ll{{\mathcal L}}
\def\vv{{\cal V}}
\def\nab{\nn^{(J,\mu)}_{ab}}
\def\mab{\mm_{ab}}
\def\ma{\mm_{ab}}
\def\lab{L_{ab}}
\def\tlab{\tilde{L}_{ab}}
\def\mabn{\mm_{ab}^N}
\def\man{\mm_a^N}
\def\labn{L_{ab}^N}
\def\hZ{\widehat{Z}}
\def\tz{\tilde{Z}}
\def\hz{\hat{z}}
\def\hg{\hat{\gamma}}
\def\mmu{m_{\mu,b,J}}
\def\tF{\tilde{F}}
\def\tf{\tilde{f}}
\def\tp{\tilde{p}}
\def\ff{{\cal F}}
\def\i{{\bf i}}
\def\jj{{\bf j}}
\def\ttau{\tilde{\tau}}
\def\tt{{\cal T}}
\def\uu{{\cal U}}
\def\hu{\widehat{U}}
\def\wloc{W_{\epsilon}}
\def\pp{{\cal P}}
\def\tg{\tilde{\gamma}}  
\def\aa{{\cal A}}
\def\tV{\widetilde{V}}
\def\cc{{\cal C}}
\def\tC{\widetilde{\cc}}
\def\Seo{S^*_\epsilon(\Omega)}
\def\sdk{S^*_{\dk}(\Omega)}
\def\lae{\Lambda_{\epsilon}}
\def\ep{\epsilon}
\def\tr{\tilde{R}}
\def\oo{{\cal O}}
\def\be{\begin{equation}}
\def\ee{\end{equation}}
\def\beqn{\begin{eqnarray*}}
\def\eeqn{\end{eqnarray*}}

\def\cc{{\mathcal C}}
\def\gi{\gamma^{(i)}}
\def\ii{{\imath }}
\def\jj{{\jmath }}
\def\tc{\tilde{C}}
\def\II{{\cal I}}
\def\ccij{ \cc_{i'_0,j'_0}[\eta]}
\def\hc{\hat{\cc}}
\def\dd{{\cal D}}
\def\la{\langle}
\def\ra{\rangle}
\def\bs{\bigskip}
\def\xio{\xi^{(0)}}
\def\xo{x^{(0)}}
\def\zo{z^{(0)}}
\def\do{\partial \Omega}
\def\dk{\partial K}
\def\dl{\partial L}
\def\cK{\hat{K}}
\def\kk{{\cal K}}
\def\pr{{\rm pr}}
\def\ff{{\cal F}}
\def\G{{\mathcal G}}
\def\dist{{\rm dist}}
\def\dds{\frac{d}{ds}}
\def\con{{\rm const}\;}
\def\Con{{\rm Const}\;}
\def\di{\displaystyle}
\def\oo{{\mathcal O}}
\def\hess{\mbox{\rm Hess}}
\def\endofproof{{\rule{6pt}{6pt}}}
\def\vm{\varphi^{(m)}}
\def\km{k^{(m)}}
\def\dm{d^{(m)}}
\def\kam{\kappa^{(m)}}
\def\dem{\delta^{(m)}}
\def\xim{\xi^{(m)}}
\def\ep{\epsilon}
\def\tt{\tilde{t}}
\def\tx{\tilde{x}}
\def\tp{\tilde{p}}
\def\eep{\hat{\epsilon}}
\def\ep{\epsilon}
\def\hu{\hat{u}}
\def\tS{\widetilde{S}}
\def\tt{{\mathcal T}}
\def\ts{\tilde{\sigma}}

\def\ms{\medskip}
\def\sn{{\SS}^{n-1}}
\def\sN{{\SS}^{N-1}}
\def\be{\begin{equation}}
\def\ee{\end{equation}}
\def\beqn{\begin{eqnarray}}
\def\eeqn{\end{eqnarray}}
\def\beqn*{\begin{eqnarray*}}
\def\eeqn*{\end{eqnarray*}}
\def\endofproof{{\rule{6pt}{6pt}}}
\def\nv{\nabla \varphi}
\def\wuloc{W^u_{\mbox{\footnotesize\rm loc}}}
\def\wsloc{W^s_{\mbox{\footnotesize\rm loc}}}

\def\dist{\mbox{\rm dist}}
\def\diam{\mbox{\rm diam}}
\def\pr{\mbox{\rm pr}}
\def\supp{\mbox{\rm supp}}
\def\Arg{\mbox{\rm Arg}}
\def\In{\mbox{\rm Int}}
\def\Im{\mbox{\rm Im}}
\def\span{\mbox{\rm span}}
\def\spec{\mbox{\rm spec}\,}
\def\Re{\mbox{\rm Re}}
\def\var{\mbox{\rm var}}
\def\conf{\mbox{\footnotesize\rm const}}
\def\Conf{\mbox{\footnotesize\rm Const}}
\def\Lip{\mbox{\rm Lip}}
\def\con{\mbox{\rm const}\;}
\def\li{\mbox{\rm li}} 
\def\ex{\mbox{\rm extd}}

\def\saa{\Sigma_A^+}
\def\sa{\Sigma_A}
\def\san{\Sigma^-_A}

\def\vxij{\varphi_{\xi,j}}
\def\tz{\tilde{z}}
\def\Xr{X^{(r)}}
\def\hs{\hat{\sigma}}
\def\hl{\hat{l}}

\def\kmax{\kappa_{\max}}
\def\kmin{\kappa_{\min}}
\def\ecc{\mbox{\rm ecc}}
\def\tB{\widetilde{B}}
\def\hh{{\mathcal H}}
\def\thh{\widetilde{\cal H}}
\def\hE{\widehat{E}}

\def\naf{\nabla f(z)}
\def\so{\sigma_0}
 \def\tp{\tilde{p}}
\def\hcij{\hat{\cc}_{i,j}}
\def\Xo{X^{(0)}}
\def\z1{z^{(1)}}
\def\Vo{V^{(0)}}
\def\Yo{Y{(0)}}
\def\tPsi{\tilde{\Psi}}
\def\hX{\hat{X}}
\def\hx{\hat{x}}

\def\pr{\mbox{\rm pr}}
\def\dK{\partial K}
\def\tpsi{\tilde{\psi}}
\def\tq{\tilde{q}}
\def\tsi{\tilde{\sigma}}
\def\iii{\sf i}
\def\clip{C^{\mbox{\footnotesize \rm Lip}}}
\def\Lip{\mbox{\rm Lip}}
\def\lip{\mbox{\footnotesize\rm Lip}}
\def\tn{t^{(n)}}
\def\En{E^{(n)}}
\def\Fn{F^{(n)}}
\def\lambdan{\lambda^{(n)}}
\def\Gr{\mbox{\rm Gr}}
\def\lan{\lambda^{(n)}}
\def\ttt{\tilde{t}}
\def\td{\tilde{d}}
\def\txi{\tilde{\xi}}
\def\tell{\tilde{\ell}}
\def\tdelta{\tilde{\delta}}
\def\tgamma{\tilde{\gamma}}
\def\tk{\tilde{k}}
\def\tu{\tilde{u}}
\def\tv{\tilde{v}}

\section{Introduction}
\renewcommand{\theequation}{\arabic{section}.\arabic{equation}}

Let $f$ be an invertible  transitive subshift of finite type over a bilateral symbol space $X$ and let $\mu$ be a Gibbs measure 
for $f$ determined by a H\"older continuous potential on $X$ (see e.g. \cite{B} or \cite{PP}). In particular, $\mu$ is $f$-invariant and ergodic.
Let $A$ be a continuous linear cocycle over $f$ acting on a continuous $\R^d$-bundle $E$ over $X$.  Thus, 
$A(x) : E(x) \longrightarrow E(f(x))$ is a linear map for each $x\in X$ and  
$$A^n(x) = A(f^{n-1}(x)) \circ A(f^{n-2}(x)) \circ \ldots \circ A(f(x)) \circ A(x) : E(x) \longrightarrow E(f^n(x)) $$
for every integer $n \geq 1$.  Let 
 $$\lambda_{k} < \lambda_{k-1} < \ldots < \lambda_2 < \lambda_1 $$
 be the {\it Lyapunov spectrum} of $f$ with respect to $\mu$, which is constant on a subset $\ll$ of $X$ of full measure.
 Denote by $E_i(x)$ the {\it Oseledets subspace} of $T_xM$ corresponding to the Lyapunov exponent $\lambda_i$.
 
Given $\ep > 0$, a measurable function  $R_\ep : \ll \longrightarrow (1,\infty)$ is called a {\it Lyapunov $\ep$-regularity function } if
$$\frac{\|v\|}{R_\ep(x)\, e^{|n|\ep}} \leq \frac{\|df^n(x)\cdot v\|}{e^{n\lambda_i}} 
\leq R_\ep(x)\, e^{|n|\ep} \|v \| \quad , \quad  v\in E_i(x) \;, \; n \in \Z ,$$
for all $x\in \ll$ and all $i = 1, \ldots, k$ and
$$e^{-\ep} \leq \frac{R_\ep (f(x))}{R_\ep (x)} \leq  e^{\ep} \quad , \quad x\in \ll .$$
A compact subset $P$ of $\ll$ is called a {\it Pesin set} if for every $\ep > 0$ there exist a  Lyapunov $\ep$-regularity function 
$R_\ep$ and a constant $C > 0$ such that $R_\ep(x) \leq C$ for all $x \in P$. A Pesin set $P$ will be called a
{\it Pesin set with exponentially small tails} for $\mu$ if for every $\delta > 0$ there exist constants $C, c > 0$ such that
$$\mu\left( \left\{ x \in \ll : \#\{ m : 0 \leq m \leq n-1\;, \; f^m(x) \notin P\} \geq \delta\, n  \right\} \right) \leq C\, e^{-c n}$$
for every integer $n \geq 0$. In other words, for every $n \geq 0$, removing a set of exponentially small measure from $\ll$,
for the remaining $x$, at least $(1-\delta)n$ points from the orbit $x, f(x), \ldots, f^{n-1}(x)$ belong to the Pesin set $P$.



\def\tF{\widetilde{F}}

Our main result in this paper is the following.
 
\bs
 
\noindent
{\sc Theorem 1.1.}
{\it Let $f : X \longrightarrow X$ be as above and let $\mu$ be a Gibbs measure on $X$  determined 
by a H\"older continuous potential on $X$. Assume that $A$ is an invertible linear cocycle over $f$ such that 
$A^{-1}$ is continuous as well and the Oseledets flags
$F_j(x) = E_j(x) \oplus E_{j-1}(x) \oplus \cdots \oplus E_1(x)$ depend upper semi-continuously on $x$. 
Then there exists a Pesin set with exponentially small tails for $\mu$.}

\bs

See Definition 2.3 in Sect. 2 below for the precise meaning of upper semi-continuity of Oseledets flags.

We will prove Theorem 1.1  in Sect. 3 below.

The motivation to study hyperbolic systems (e.g. Axiom flows on basic sets) admitting Pesin sets with exponentially small tails
comes from \cite{St} which deals with the case of contact Anosov flows. It is proved there that every Gibbs measure 
for such a flow admitting a Pesin set with exponentially small tails has exponential mixing.

\def\hF{\widehat{F}}

\section{Lyapunov exponents and Oseledets subspaces}
\renewcommand{\theequation}{\arabic{section}.\arabic{equation}}
\setcounter{equation}{0}

To prove Theorem 1.1 we will use certain background related to the so called Multiplicative
Ergodic Theorem  and some arguments from \cite{GS}. 

 
 Let $f$ be an invertible  transitive subshift of finite type over a bilateral symbol space $X$ and let $\mu$ be a Gibbs measure 
for $f$ determined by a H\"older continuous potential on $X$ (see e.g. \cite{B} or \cite{PP}). In particular, $\mu$ is $f$-invariant and ergodic.
Let $A$ be a continuous linear cocycle over $f$ acting on a continuous $\R^d$-bundle $E$ over $X$.  Thus, 
$A(x) : E(x) \longrightarrow E(f(x))$ is a linear map for each $x\in X$ and  
$$A^n(x) = A(f^{n-1}(x)) \circ A(f^{n-2}(x)) \circ \ldots \circ A(f(x)) \circ A(x) : E(x) \longrightarrow E(f^n(x)) $$
for every integer $n \geq 1$.

Given an integer $p = 1, \ldots,d$, let $\Gr_p(\R^d)$ be the {\it Grassman manifold}  of the linear subspaces  of $\R^d$ of dimension $p$
endowed with the usual distance $d(U,V)$ between subspaces $U,V \in \Gr_p(\R^d)$ defined by
$$d(U,V) =  \max \{ |\la u, w\ra | : u \in U, w\in V^\perp , \|u\| = \|w\| = 1\} . $$

The following is Oseledets' Multiplicative Ergodic Theorem stated under the above assumptions (see  e.g.  
\cite{Ar}, \cite{BP}, \cite{V} or \cite{Sa}  for  related detailed exposition and proofs).

\bs

\noindent
{\sc Theorem 2.1.}
(Multiplicative Ergodic Theorem). 
{\it There exists a subset $\ll$ of $X$ with  $\mu(\ll) = 1$ such that:}

\ms

(a) {\it  For all $x\in \ll$ there exists}
$\displaystyle N(x) = \lim_{n\to\infty} (A^n(x)^t A^n(x))^{1/2n} .$

\ms

(b) {\it There exist an integer $k \geq 1$  
such that  the operator (matrix) $N(x)$ has $k$ distinct eigenvalues $t_k(x) < \ldots < t_{1}(x)$ for all $x \in \ll$.  }

\ms

(c) {\it There exist numbers $t_k < t_{k-1} < \ldots < t_2 <  t_{1}$ such that  
$(\tn_i(x))^{1/n} \to t_i$ for all $x \in \ll$ and all $i = 1, \ldots, k$.}

\ms

(d) {\it For all $x\in \ll$ and every $j = 1, \ldots, k$ the dimension $\dim(\En_j(x)) = m_j$ is constant
and there exists $\lim_{n\to \infty} \En_j(x) = E_j(x)$ in $\Gr_{m_j}(\R^{d})$.} 

\bs

(e) {\it For all $x\in \ll$ and every $j = 1, \ldots, k$ we have
$$\lim_{n\to \infty} \frac{1}{n} \log \|A^n(x) \cdot u\| = \log t_j \quad , \quad u \in \tF_j(x)\setminus \tF_{j+1} (x) ,$$
where $\tF_{k+1}(x)  = \{ 0\}$ and $\tF_j(x) = E_k(x) \oplus E_{k-1}(x) \oplus \ldots \oplus E_{j+1}(x) \oplus E_j(x)$
for all $j = 1, \ldots, k$.}

\bs

The numbers 
$$\lambda_k = \log t_k < \lambda_{k-1} = \log t_{k-1} < \ldots  < \lambda_2 = \log t_2 < \lambda_{1} = \log t_{1}$$ 
are the   (distinct) {\it Lyapunov exponents} of the cocycle $A$ over $f$. 
Setting 
$$\lan_i(x) = \frac{1}{n} \log \tn_i(x) ,$$
we have $\lan_i(x) \to \lambda_i$ for all $x \in \ll$.  The so called
{\it Oseledets subspaces} $E_i(x)$ form the corresponding {\it Oseledets flags}
$$F_j(x) = E_j(x) \oplus E_{j-1}(x) \oplus \ldots  \oplus E_1(x) ,$$
$x \in \ll$, $j = 1, \ldots, k$. Both are invariant with respect to $A$.

\bs

Given a finite-dimensional real inner product space $V$ with inner product $\la \cdot, \cdot \ra$, consider the 
{\it $r$-fold exterior product} $\Lambda^r V$ which is naturally identified with linear combinations of wedge products
$v_1\wedge v_2 \wedge \cdots\wedge v_r$ of vectors in $V$ (see e.g. Sect. 3.2.3 in \cite{Ar}).
A natural inner product is then defined on $\Lambda^rV$ so that
$$\la u_1 \wedge \cdots \wedge u_r , v_1 \wedge \cdots \wedge v_r \ra = \det (\la u_i, v_j\ra)_{i,j=1}^r .$$
Thus, $\| u_1 \wedge \cdots \wedge u_r\|^2 = \det (\la u_i, u_j\ra)_{i,j=1}^r$.
If $L : V \longrightarrow V$ is a linear operator, then the linear operator $\Lambda^r L : \Lambda^r V \longrightarrow \Lambda^r V$ 
is defined by
$$(\Lambda^r L)  (v_1 \wedge v_2 \wedge \cdots \wedge v_r) = Lv_1 \wedge Lv_2 \wedge \cdots\wedge Lv_r.$$

Let $\chi_1 \geq \chi_2 \geq \ldots \geq \chi_{d}$ be the {\it Lyapunov exponents} $\lambda_1, \lambda_2, \ldots, \lambda_{k}$ 
of the linear cocycle $A(x)$ with the required number of repetitions,
so that for every $i  = 1,\ldots, d$ we have $\chi_i = \lambda_j$ for some (unique) $j = j_i = 1, \ldots, k$.
It is well-known that
$$\log \|\Lambda^i A(x)\| = \chi_1+  \chi_{2} +  \ldots + \chi_{i} $$
for all $i = 1, \ldots, d$ (see e.g. Proposition 3.2.7 in \cite{Ar}).

\bs

\noindent
{\sc Definition 2.2.}\label{Def1}
(\cite{GS})
We say that the cocycle $A$ {\it has exponential large deviations for all exponents} if for
any $\epsilon>0$ there exist constants $C, c > 0$  such that
$$  \mu \left(\left\{  x \in X :  \left| \log \|\Lambda^i A^n(x)\| - n (\chi_1 + \ldots + \chi_{i}) \right| \geq n \epsilon\right \}\right) \leq C e^{-c n}$$
 for all $n\geq 0$ and all $i  = 1, \ldots, d$.

\bs

It follows from Proposition 3.3 in \cite{GS}  that  
\begin{equation}\label{eq:2.1}
\mu \left(\left\{ x \in \ll : \lambdan_1(x) + \ldots + \lambdan_i(x) \geq \lambda_1 + \ldots + \lambda_i + \ep \right\} \right)\leq C e^{-cn}
\end{equation}
for all $n \geq 1$. That is, given $i = 1,\ldots, k$, for every $n \geq 1$ we have 
$$\lambdan_1(x) + \ldots + \lambdan_i(x) < \lambda_1 + \ldots + \lambda_i + \ep$$
for `most' $x$.

We will prove a similar estimate from above for every Lyapunov exponent $\lambda_i$ under 
a certain additional assumption. 

\bs
 
\noindent
{\sc Definition 2.3.}\label{Def2}
We will say that the {\it Oseledets flags $F_j(x)$ related to $\mu$ depend upper semi-continuously on $x$} if we can choose the set $\ll$ 
in Theorem 2.1 with $\mu(\ll) = 1$ so that for every $x \in \ll$, every sequence $\{x_n\}$ in $\ll$ converging to $x$, every $j = 1, \ldots, k$ and every 
sequence $\{u_n\}_{n=1}^\infty$ with $u_n \in F_j(x_n)$  for all $n \geq 1$, if there exists $u = \lim_{n\to\infty} u_n$,  then $u \in F_j(x)$.

\bs

In the next theorem the assumption that the Oseledets flags $F_j(x)$ depend upper semi-continuously on $x$ is essential.
As the counterexamples in section A in \cite{GS} show, such a statement cannot be true without additional assumptions.


\bs
 
\noindent
{\sc Theorem 2.4.}
 {\it Let $(f, \mu)$ be an invertible  transitive subshift of finite type with a Gibbs measure $\mu$ and let $A$ be an invertible
continuous linear cocycle above $f$ with Lyapunov exponents \\ $\lambda_k < \lambda_{k-1}  \ldots < \lambda_2 <  \lambda_1$. 
Assume that $A^{-1}$ is also continious and the Oseledets flags $F_j(x)$ depend upper semi-continuously on $x$. Then for any $\ep > 0$   
there exist constants $C, c > 0$ such that 
\begin{equation}\label{eq:2.2}
\mu \left(\left\{ x \in \ll : \lan_j(x)  < \lambda_j -  \ep \right\} \right) \leq C e^{-cn} 
\end{equation}
for all $n \geq 0$ and all $j = 1, \ldots, k$.}

\bs

\noindent
{\sc Corollary 2.5.}
{\it Under the assumptions of Theorem 2.4,  for every $\ep > 0$ there exist constants $C, c > 0$ such that 
\begin{equation}\label{eq:2.3}
\mu \left(\left\{ x \in \ll : \left|\lan_j(x) - \lambda_j \right| \geq  \ep \right\} \right) \leq C e^{-cn} 
\end{equation}
for all $n \geq 0$ and all $j = 1, \ldots, k$.}

\bs

\noindent
{\it Proof of Corollary } 2.5. Given $\ep > 0$, take $C, c > 0$ so that (\ref{eq:2.1})  and (\ref{eq:2.2}) hold with
$\ep$ replaced by $\ep/k$ and $C$ replaced by $C/(2k)$. For every $j = 1, \ldots, k$ set
$$X_j = \left\{ x \in \ll : \left|\lan_j(x) - \lambda_j \right| \geq  \ep \right\} ,$$
$$Y_j = \left\{ x \in \ll : \lan_j(x) < \lambda_j  - \ep/k \right\} ,$$
$$Z_j = \left\{ x \in \ll : \lambdan_1(x) + \ldots + \lambdan_j(x) \geq \lambda_1 + \ldots + \lambda_j + \ep/k \right\} .$$
Then $\mu(Y_j) \leq C e^{-cn}/(2k)$ and $\mu(Z_j) \leq C e^{-cn}/(2k)$ for all $j = 1, \ldots, k$ and all $n \geq 0$.
On the other hand,
\begin{equation}\label{eq:2.4}
X_j \subset Z_j \cup \cup_{i=1}^k Y_i
\end{equation}
for all $j$. Indeed, given $j$ and $x \in X_j$, assume that $x \notin \cup_{i=1}^k Y_i$. This  gives
\begin{equation}\label{eq:2.5}
\lan_i(x) \geq \lambda_i  - \ep/k \quad, \quad i = 1, \ldots, k ,
\end{equation} 
while $x \in X_j$ implies $\lan_j(x) \geq \lambda_j + \ep$ 
or $\lan_j(x) \leq \lambda_j - \ep$. The latter is impossible by $x \notin Y_j$, so we must have
$\lan_j(x) \geq \lambda_j + \ep$. Using this and (\ref{eq:2.5}) for $i = 1,\ldots,j-1$ we get
\begin{eqnarray*}
&         &\lan_1(x) + \lan_2(x) + \ldots + \lan_j(x)\\
& \geq & (\lambda_1 -\ep/k) + (\lambda_2 - \ep/k) + \ldots + (\lambda_{j-1} - \ep/k) + \lan_j(x)\\
& \geq & \lambda_1 + \lambda_2 + \ldots + \lambda_{j-1} - (j-1)\ep/k + (\lambda_j + \ep)\\
& >      & \lambda_1 + \lambda_2 + \ldots + \lambda_{j-1} + \lambda_j + \ep/k .
\end{eqnarray*}
Thus, $x \in Z_j$. This proves (\ref{eq:2.4}), which immediately implies $\mu(X_j) \leq C e^{-cn}$ for all $n \geq 0$.
\endofproof

\bs

\noindent
{\it Proof of Theorem 1.1}. 
It follows from Corollary 2.5 that the cocycle $A$  has exponential large deviations for all
exponents.  Combining this with Theorem 1.7 in \cite{GS} yields Theorem 1.1.
\endofproof



\section{Proof of Theorem 2.4}
\setcounter{equation}{0}

 Let again $\mu$ be a Gibbs measure on $X$ defined by a H\"older continuous potential and let 
 $$\lambda_k <  \lambda_{k-1} < \ldots < \lambda_2 < \lambda_1$$
be the Lyapunov spectrum of the cocycle $A$ with respect to $\mu$, constant on a measurable subset $\ll$ of $X$ of full measure.
Consider the corresponding invariant decomposition
$$E_k(x) \oplus E_{k-1}(x)  \oplus \cdots \oplus E_2(x) \oplus E_{1}(x)  = \R^d  ,$$
where $E_i(x)$ is the Oseledets subspace corresponding to the Lyapunov exponent $\lambda_i$. 
We assume that the Oseledets flags 
$$F_j(x) = E_j(x) \oplus E_{j-1}(x) \oplus \ldots  \oplus E_1(x)$$
depend  upper semi-continuously on $x \in \ll$ (see Definition 2.3).



To derive (\ref{eq:2.2}) we will use the the cocycle $B = A^{-1}$ over the inverse transformation \\
$f^{-1} : X \longrightarrow X$. $B$ is continuous by assumption and its Lyapunov exponents are 
$$-\lambda_1 <  - \lambda_{2}  <  \ldots < - \lambda_{k-1} <  -\lambda_k $$
with the same Oseledets subspaces $E_i(x)$  (see e.g. Sect. 2.6.3 in \cite{Sa}).

{\bf Fix $j = 1, \ldots, k$. Throughout this section $j$ will stay fixed.}
Set
$$\gamma_n(x) = \log \| B^n(x)_{| F_j(x)}\| \quad , \quad x\in \ll \:, \: n \geq 1 .$$

{\bf Fix an  $\ep > 0$.} We will now prove that  there exist constants $C, c > 0$ such that 
\begin{equation}\label{eq:3.1}
\mu \left(\left\{ x \in \ll : \frac{1}{n} \gamma_n (x)  > - \lambda_j + \ep \right\} \right) \leq C e^{-cn} 
\end{equation}
for all $n \geq 0$. From this (\ref{eq:2.2}) follows immediately.

Clearly, $\gamma_n(x) = \log \| B^n(x)_{| E_j(x)}\|$, $\gamma_n$ is a sub-additive function on $\ll$, and
$|\gamma_n(x)| \leq \Con < \infty$ for all $x\in \ll$.   We can assume $\ll$ is so that 
$$\lim_{n\to\infty} \frac{\gamma_n(x)}{n} = - \lambda_j$$
for all $x \in \ll$.  By Kingsman's Ergodic Theorem (see e.g. \cite{Ar}, \cite{V} or \cite{Sa}) there exists the limit
\be\label{eq:3.2}
\lim_{n\to \infty} \frac{1}{n} \int_X \gamma_n(x) \, d\mu(x) =
\inf_{n} \frac{1}{n} \int_X \gamma_n(x) \, d\mu(x) = -\lambda_j .
\ee
Set
$$\hF_j(x) = \{ u \in F_j(x) : \|u\| = 1 \} .$$
Since
$$\gamma_n(x) =  \max \left\{ \log \left\| B^n(x) \cdot v\right\| : v\in \hF_j(x) \right\}  , $$
we will derive from our assumptions that $\gamma_n$ is upper semi-continuous on $\ll$.
However it is not clear whether $\gamma_n$ is continuous and whether it has a continuous extension to $X$, and so it is not immediately
possible to use Lemma 3.2\footnote{Lemma 3.2 (\cite{GS}): {\it Let $(T, \mu)$ be a transitive subshift of finite type with a Gibbs measure $\mu$
and let $a(n,x)$ be a subadditive cocycle above $T$ such that $a(n,x)$ is continuous for all $n$. Let $\chi$ be
the almost sure limit of $a(n,x)/n$ and assume that $\chi > -\infty$. Then, for any $\epsilon > 0$, there
exist constants $C > 0$, $c > 0$ such that, for all $n \geq 0$ we have}
$\mu \left(\{ x\in X:  a(n,x)/n \geq  \chi + \epsilon\}\right) \leq C\, e^{-c n} .$} in \cite{GS}.

\bs
 
\noindent
{\sc Lemma 3.1.}
{\it For all $n \geq 1$, the function $\gamma_n$ is upper semi-continuous on $\ll$, i.e.
$$U = \{ x \in \ll : \gamma_n(x) < \omega\}$$
is an open subset of $\ll$ for all $\omega \in \R$.}

\bs

\noindent
{\it Proof.} Let $n \geq 1$ and let  $U$ be defined as above for some $\omega \in \R$.  We will prove that $\ll \setminus U$ is closed in $\ll$.
Let $x_m \to x$ in $\ll$ for some sequence $\{x_m\} \subset \ll \setminus U$. Then
$\gamma_n(x_m) \geq \omega$ for all $m$. From the definition of $\gamma_n$, 
for all $m \geq 1$ there exists $u_m \in \hF_j(x_m)$ with $\gamma_n(x_m) = \log \|B^n(x_m)\cdot u_m\|$.
Replacing $\{x_m\}$ with an appropriate subsequence, we may assume that $u_m \to u$ for some $u \in \R^d$.
Then $\|u \| = 1$ and the upper semi-continuity of the bundle $F_j$ (see Definition 2.3) implies that $u \in F_j(x)$.
We now get
$$\omega \leq \gamma_n(x_m) = \log \|B^n(x_m)\cdot u_m\| \to \log \|B^n(x)\cdot u\| \leq \gamma_n(x) .$$
Therefore $\gamma_n(x) \geq \omega$, so $x \in \ll \setminus U$. Thus, $\ll \setminus U$
is closed in $\ll$ and so $U$ is open in $\ll$.
\endofproof

\bs

\noindent
{\it Proof of Theorem} 2.4. Since $B$ is continuous, there exists a constant $C_0 > 0$ such that $\|B(x)\| \leq C_0$ for all $x \in X$.
Then $ 0 \leq \gamma_n(x) \leq n C_0$ for all $x\in \ll$ and all $n \geq 1$.

Since $\gamma_n$ is bounded it follows e.g. from Theorem 2.2 in \cite{Kos} that it has an 
upper semi-continuous  extension to $X$ (which again takes values in a bounded interval).
We assume below that for every $n \geq 1$, $\tgamma_n$ is such an upper semi-continuous  extension of the
original function $\gamma_n$ to $X$.

By (\ref{eq:3.2}), there exists $n_0 \geq 1$ so that
$$\int_X \frac{\tgamma_n(x)}{n}\, d\mu(x) = \int_X \frac{\gamma_n(x)}{n}\, d\mu(x) < - \lambda_j + \frac{\ep}{8} $$
for all $n \geq n_0$. 

{\bf Fix an $m \geq n_0$.}
It follows from Baire's Theorem for semi-continuous functions (\cite{Ba}) and the upper semi-continuity
of $\tgamma_m$ that there exists a sequence 
$$h_1 \geq h_2 \geq \ldots \geq h_n \geq \ldots$$
of continuous functions on $X$ such that $h_n(x) \to \tgamma_m(x)$ as $n \to \infty$ for all $x \in X$.
Now Lebesgue's Dominated Convergence Theorem implies
$$\lim_{n\to \infty} \int_X \frac{h_n}{m}\, d\mu  = \int_X \frac{\tgamma_m}{m}\, d\mu 
= \int_X \frac{\gamma_m}{m}\, d\mu < - \lambda_j +  \frac{\ep}{8},$$
therefore for sufficiently large $n$ we have
$$\int_X \frac{h_n}{m}\, d\mu  <  - \lambda_j + \frac{\ep}{8} .$$
{\bf Fix such an $n \geq 1$.} It follows from well-known results about approximation of continuous functions by
Lipschitz ones (see e.g. \cite{Ge}) that there exists a Lipschitz function $H$ on $X$ so that 
$$|H(x) - h_n(x)| < \frac{\ep}{8}$$
for all $x\in X$. Set $G_m = H + \ep/8$. Then $G_m$ is Lipschitz and
$$G_m \geq (h_n - \ep/8)  + \ep/8 = h_n \geq \tgamma_m .$$
In particular,
\begin{equation}\label{eq:3.3}
G_m(x) \geq \gamma_m(x) \quad, \quad x \in \ll .
\end{equation}
Moreover
\begin{equation}\label{eq:3.4}
\int_X \frac{G_m}{m}\, d\mu = \int_X \frac{H + \ep/8}{m}\, d\mu \leq  \int_X \frac{h_n + \ep/4}{m}\, d\mu  < - \lambda_j + \frac{3\ep}{8} ,
\end{equation}
assuming $m \geq 2$.

\bs

Given arbitrary integers $n, s \geq 1$ and $x \in \ll$, let $u\in \hF_j(x)$ be such that 
$$\gamma_{n+s}(x) = \log \left\| B^{n+s}(x) \cdot u\right\|, $$
that is 
$$\log \left\| B^{n+s}(x) \cdot u\right\| = \max \left \{ \log \left\| B^{n+s}(x) \cdot v\right\| : v\in \hF_j(x) \right\} .$$
Set
$$w = \frac{1}{\left\| B^s(x) \cdot u\right\|} \, B^s(x) \cdot u \in \hF_j(f^s(x)) .$$ 
We now derive
\begin{eqnarray*}
\gamma_{n+s}(x)
& =     & \log \left\| B^{n+s}(x) \cdot u\right\| = \log \left\| B^n(f^s(x)) \cdot ( B^s(x) \cdot u)|\right\|\\
& =     & \log \left(\left\| B^n(f^s(x)) \cdot w\right\| \cdot  \|B^s(x) \cdot u \|\right)\\
& =     & \log \left\| B^n(f^s(x)) \cdot w\right\| + \log  \|B^s(x) \cdot u \| \\
& \leq &  \gamma_n(f^s(x)) +  \gamma_s(x) .
\end{eqnarray*}
Thus
\begin{equation}\label{eq:3.5}
\gamma_{n+s}(x) \leq  \gamma_n(f^s(x)) + \gamma_s(x) 
\end{equation}
for all $x \in \ll$ and all integers $n,s \geq 1$.

Recall the {\bf fixed integer} $m \geq n_0$.
Clearly $\gamma_s(x) \leq D = m C_0$ for all $x \in \ll$ and all $s = 0,1,\ldots, m$.

Next, we will use  arguments similar to the ones in the proofs of Lemmas 3.1 and 3.2 in \cite{GS}.

Let $x \in \ll$ and let $n \geq 1$. Fix for a moment an arbitrary integer $\ell$ with $0 \leq \ell \leq m-1$.
Then $n = r+ mp + \ell$ for some integers $p$ and $r$ with $0 \leq r \leq m-1$. Applying (\ref{eq:3.5}) inductively 
and using (\ref{eq:3.3}), we get
\begin{eqnarray*}
\gamma_n(x)
& =     & \gamma_{(r+mp) + \ell}(x) \leq \gamma_{r+mp} (f^\ell(x)) + \gamma_\ell(x)
 \leq  \gamma_{(r+ m(p-1))+ m}(f^\ell(x)) + D\\
& \leq & \gamma_{r+ m(p-1)}(f^{m+\ell}(x)) + \gamma_m(f^\ell(x)) + D\\
& \leq &  \gamma_{r+ m(p-1)}(f^{m+\ell}(x)) + G_m(f^\ell(x)) + D\\
& \leq & \gamma_{r+ m(p-2)}(f^{2m+\ell}(x)) + G_m(f^{m+\ell}(x)) + G_m(f^\ell(x)) + D \leq \ldots \\
& \leq & \gamma_{r}(f^{mp+\ell}(x)) + \sum_{q=0}^{p-1} G_m(f^{mq + \ell}(x)) + D
\leq \sum_{q=0}^{p-1} G_m(f^{mq + \ell}(x)) + 2 D .
\end{eqnarray*}
Hence
$$\gamma_n(x) \leq  \sum_{q=0}^{p-1} G_m(f^{mq + \ell}(x)) + 2 D .$$
Summing up over $\ell = 0,1, \ldots, m-1$ this gives
$$m\, \gamma_n(x) \leq  \sum_{\ell=0}^{m-1}\sum_{q=0}^{p-1} G_m(f^{mq + \ell}(x)) + 2 m D 
\leq (G_m)_n(x) + 3m D .$$
Thus, assuming that the integer $n_0 \geq 1$  is chosen sufficiently large, we have
\begin{equation}\label{eq:3.6}
\frac{1}{n} \gamma_n(x) \leq \frac{1}{n}\, \left( \frac{G_m}{m} \right)_n(x) + \frac{3D}{n} < \frac{1}{n}\, \left( \frac{G_m}{m} \right)_n(x) + \frac{\ep}{4} 
\end{equation}
for all $x\in \ll$ and all $n \geq n_0$.

Set
$$\Gamma_n = \left\{ x \in \ll : \frac{1}{n} \gamma_n(x)  > - \lambda_j + \ep \right\} ,$$
$$Y_n = \left\{ x \in \ll : \frac{1}{n}  \left( \frac{G_m}{m} \right)_n(x)  \geq \int_X \frac{G_m}{m}\, d\mu + \frac{\ep}{4} \right\} .$$
Since $G_m$ is Lipschitz, we can now use the classical {\bf Large Deviation Principle} (see \cite{Y} or \cite{Ki})
for H\"older continuous functions on $X$. It implies that there exist constants $C, c > 0$ such that
\begin{equation}\label{eq:3.7}
\mu(Y_n) \leq C\, e^{-c n} \quad, \quad n \geq 1 .
\end{equation}
On the other hand, $\Gamma_n \subset Y_n$ for $n \geq n_0$. Indeed, given $n \geq n_0$ and $x \in \Gamma_n$, 
(\ref{eq:3.6}) and (\ref{eq:3.4}) imply
\begin{eqnarray*}
\frac{1}{n}\, \left( \frac{G_m}{m} \right)_n(x)  
& \geq &\frac{1}{n} \gamma_n(x) - \frac{\ep}{4} \geq (-\lambda_j + \ep) - \frac{\ep}{4} = - \lambda_j + \frac{3\ep}{4}\\
&  >     & \left(\int_X \frac{G_m}{m}\, d\mu - \frac{3\ep}{8}\right) + \frac{3\ep}{4}  >   \int_X \frac{G_m}{m}\, d\mu + \frac{\ep}{4} ,
\end{eqnarray*}
so $x \in Y_n$. Now (\ref{eq:3.7}) yields $\mu(\Gamma_n) \leq C e^{-cn}$ for all $n \geq n_0$. Taking a larger constant $C$ if necessary,
we get $\mu(\Gamma_n) \leq C e^{-cn}$ for all $n \geq 1$. 
This proves (\ref{eq:3.1}). As we mentioned already, (\ref{eq:2.2}) follows immediately.
This concludes the proof of Theorem 2.4.
\endofproof

\def\beqn{\begin{eqnarray*}}
\def\eeqn{\end{eqnarray*}}
\def\yo{y^{(0)}}
\def\tmu{\tilde{\mu}}

\bs

\noindent
{\it Acknowledgments.} Thanks are due to Sebastien Gou\"ezel for his comments on earlier drafts of this paper.

{\footnotesize


\begin{thebibliography}{4}


\bibitem[Ar]{Ar}
\newblock L. Arnold.
\newblock \emph{Random Dynamical Systems}.
\newblock Springer Monographs in Mathematics, Springer, Berlin 1998.


\bibitem[Ba]{Ba}
\newblock R. Baire.
\newblock \emph{Lecons sur les fonctions discontinues, profess\'ees au coll\'ege de France}.
\newblock Gauthier-Villars (1905).

\bibitem[BP]{BP}  
\newblock L. Barreira and Ya. Pesin.
\newblock \emph{Lyapunov exponents and smooth ergodic theory}.
\newblock Univ. Lect. Series 23, American Mathematical Society,  Providence,  RI, 2001.


\bibitem[B]{B}  
\newblock R. Bowen.
\newblock \emph{Equilibrium states and the ergodic theory of Anosov diffeomorphisms}.
\newblock Lect. Notes in Maths. \textbf {470}, Springer-Verlag, Berlin, 1975.











\bibitem[Ge]{Ge}  
\newblock  G. Georganopoulos.
\newblock {Sur l'approximation des fonctions continues par des fonctions lipschitziennes}.
\newblock \emph{C. R. Acad. Sci. Paris S\'er. A-B}  \textbf{264}  (1967), 319-321.

\bibitem[GS]{GS}  
\newblock S. Gou\"ezel and L. Stoyanov.
\newblock {Quantitative Pesin theory for Anosov diffeomorphisms and flows}.
\newblock \emph{ Ergodic Th. \& Dyn. Sys.} \textbf{39} (2019), 159-200.

\bibitem[Ha]{Ha}  
\newblock B. Hasselblatt.
\newblock {Regularity of the Anosov splitting and of horospheric foliations}.
\newblock \emph{Ergod. Th.\& Dynam. Sys.} \textbf{14} (1994), 645-666.





\bibitem[KH]{KH}  
\newblock A. Katok and  B. Hasselblatt.
\newblock \emph{Introduction to the Modern Theory of Dynamical System}.
\newblock Cambridge Univ. Press, Cambridge 1995.


\bibitem[Ki]{Ki}  
\newblock Yu. Kifer.
\newblock {Large deviations in dynamical systems and stochastic processes}.
\newblock \emph{Trans. Amer. Math. Soc. } \textbf{321} (1990), 505-524.




\bibitem[Kos]{Kos}  
\newblock J. Kosman.
\newblock {Extensions of semicontinuous and quasicontinuous functions from dense subspaces}.
\newblock \emph{Quaestiones Math.} \textbf{43} (2020), 1385-1390.






\bibitem[PP]{PP}  
\newblock W. Parry and M. Pollicott.
\newblock {Zeta functions and  the periodic orbit structure of hyperbolic dynamics}.
\newblock \emph{Ast\'erisque} \textbf{187-188}, (1990).


\bibitem[PS]{PS}  
\newblock V. Petkov and L. Stoyanov.
\newblock \emph{Geometry of the generalized geodesic flow and inverse spectral problems}.
\newblock 2nd edition. John Wiley \& Sons, Chichester, 2017.

\bibitem[Sa]{Sa}  
\newblock O. Sarig.
\newblock \emph{Lecture notes on ergodic theory}.
\newblock Lecture Notes, Penn. State University, 2009.

\bibitem[St]{St}  
\newblock L. Stoyanov.
\newblock {Spectral properties of Ruelle transfer operators for regular Gibbs 
measures  and decay of correlations for  contact Anosov flows}.
\newblock \emph{Memoirs Amer. Math. Soc. } vol. \textbf{283}, No. 1404 (2023); arXiv:1712.03103.

\bibitem[V]{V}  
\newblock M. Viana.
\newblock \emph{Lectures on Lyapunov exponents}.
\newblock Cambridge Studies in Adv. Math. vol.145, Cambridge Univ. Press 2014.  



\bibitem[Y]{Y}  
\newblock L.-S. Young.
\newblock {Large deviations in dynamical systems}.
\newblock \emph{Trans. Amer. Math. Soc.} \textbf{318} (1990), 525-543.




\end{thebibliography}
\end{document}